# Application of Differential Transform Method for El Nino Southern Oscillation (ENSO) Model with compared Adomian Decomposition and Variational Iteration Methods

Murat Gubes[a], H. Alpaslan Peker[b], Galip Oturanc[b]

[a]Department of Mathematics, Kamil Ozdag Science Faculty of Karamanoglu Mehmetbey University, 70100, Campus, Karaman, Turkey

[b]Department of Mathematics, Science Faculty of Selcuk University, 42100, Campus-Konya, Turkey

E-mail: mgubes@kmu.edu.tr, hapeker@selcuk.edu.tr, goturanc@selcuk.edu.tr

*Abstract*

We consider two nonlinear El Nino Southern Oscillation (ENSO) model to obtain approximate solutions with differential transform method for the first time. Efficiency, accuracy and error rates of solutions are compared with analytic solution, variational iteration and adomian decomposition solutions on the given tables and figures.

## 1. Introduction

In the past fifty years, nonlinear problems which appeared physical phenomena, engineering applications and various of scientific areas are modelled and they are investigated by using so many computable algorithms and approximating methods. Some of these numerical methods are Homotopy Perturbation Method (HPM), Homotopy Analysis Method (HAM), Adomian Decomposition Method (ADM), Variational Iteration Method (VIM) and Differential Transformation Method (DTM). Many authors studied nonlinear models to compute approximate solutions and their convergences with these methods (Ref. [1]-[7] ).

One of the nonlinear mathematical models is El Nino/La Nina Southern Oscillation (ENSO) model which is a climatic pattern that appear in the atmosphere and ocean of tropical areas. It occurs with two phases, the warm oceanic phase El Nino and the cold phase La Nina. The Natural ENSO event has a prolonged effect on the global climate. Therefore, the nonlinear ENSO models are very important and interesting subject to investigate ocean climate, atmospheric physics and dynamical systems (see more [8]-[21]).

The Enso delayed oscillator, delay sea-air oscillator model, atmosphere-ocean model and sea-air coupled dynamical system were solved approximately, simulated graphically with perturbation, ADM, Modified VIM, Homotopy, Asymtotic and Laplace-Adomian-Pade Techniques by authors in Ref. [8]-[21].

In this paper, we consider the two ENSO Models and compute the approximate solution and error rates by using differential transform method (DTM). Unlike the adomian decomposition method (ADM), variational iteration method (VIM) and homotopy perturbation method (HPM), the DTM transforms the nonlinear models into algebraic equations. Therefore, DTM is applied linear and nonlinear differential equations very easily. So, we perform the DTM for two ENSO models.

Numerical results related to these models are shown to verify the efficiency and accuracy of the DTM compared with the ADM and VIM solutions.

## 2. Basic Definitions of Differential Transform Method (Dtm)

Differential transform method is a numerical method based on Taylor expansion. This method tries to find coefficients of series expansion of unknown function by using the initial data on the problem. The concept of differential transform method was first proposed by Zhou (Ref. [2]).

**Definition 1:** The one-dimensional differential transform of a function $y(x)$ at the point $x=x_0$ is defined [2],[7]

$$Y(k) = \frac{1}{k!}[\frac{d^k}{dx^k} y(x)]_{x=x_0} \tag{1}$$

where $y(x)$ is the original function and $Y(k)$ is the transformed function.

**Definition 2:** The differential inverse transform of $Y(k)$ is defined as follows (see Ref. [2],[7])

$$y(x) = \sum_{k=0}^{\infty} Y(k)(x-x_0)^k \tag{2}$$

From (1) and (2), we write down

$$y(x) = \sum_{k=0}^{\infty} \frac{1}{k!}[\frac{d^k}{dt^k} y(x)]_{x=x_0} (x-x_0)^k \tag{3}$$

Related the above definitions and from Ref. [2] and [7], we present some basic theorems of DTM as follow:

**Theorem 1:** If $f(x) = g(x) \pm h(x)$, then transformation form $F(k) = G(k) \pm H(k)$.

**Theorem 2:** If $f(x) = cg(x)$, then transformation form $F(k) = cG(k)$.

**Theorem 3:** If $f(x) = \frac{d^n g(x)}{dx^n}$, then transformation form $F(k) = \frac{(k+n)!}{(k)!} G(k+n)$.

**Theorem 4:** If $f(x) = g(x)h(x)$, then transformation form

$$F(k) = G(k) \otimes H(k) = \sum_{r=0}^{k} G(r)H(k-r).$$

**Theorem 5:** If $f(x) = g(x)h(x)v(x)$, then transformation form

$$F(k) = G(k) \otimes H(k) \otimes V(k) = \sum_{r=0}^{k}\sum_{t=0}^{k-r} G(r)H(t)V(k-r+t).$$

**Theorem 6:** If $f(x) = g(x)\dfrac{d^2h(x)}{dx^2}$, then transformation form

$$F(k) = G(k) \otimes \left[(k-r+1)(k-r+2)H(k+2)\right]$$
$$= \sum_{r=0}^{k}(k-r+1)(k-r+2)G(r)H(k-r+2)$$

**Theorem 7:** If $f(x) = ax^m$, $a \in \mathbb{R}$, then transformation form

$$F(k) = a\delta(k-m) = \begin{cases} a, & k = m \\ 0, & otherwise \end{cases}.$$

The proofs of above theorems are shown in Ref. [2] and [7].

**Definition 3**: When there exists nonlinear terms with respect to nonlinear functions, then the differential transform of these terms are defined

$$N(k) = \frac{1}{k!}\left[\frac{d^k}{dx^k}Ny(x)\right]_{x=x_0} \tag{4}$$

Here, $Ny(x)$ is the main nonlinear term and $N(k)$ is transformed form.

**Definition 4:** With similar way from (2), differential inverse transform of $N(k)$ is written as

$$Ny(x) = \sum_{k=0}^{\infty} N(k)(x-x_0)^k \tag{5}$$

Thus, combining (4) and (5) we have

$$\begin{aligned} Ny(x) &= \sum_{k=0}^{\infty}\frac{1}{k!}\left[\frac{d^k}{dx^k}Ny(x)\right]_{x=x_0}(x-x_0)^k \\ &= \sum_{k=0}^{\infty}\frac{1}{k!}\left[\frac{d^k}{dx^k}N\left(\sum_{r=0}^{\infty}Y(r)(x-x_0)^r\right)\right]_{x=x_0}(x-x_0)^k \end{aligned} \tag{6}$$

In real applications, $Ny(x)$ is expressed as a finite series

$$Ny(x) = \sum_{k=0}^{m} N(k)(x-x_0)^k \tag{7}$$

## 3. Enso Model And Approximate Solution By Dtm

**Problem 1:** Firstly, the coupled dynamical ENSO system was considered to describe the oscillating physical mechanism as follows (see Ref. [8],[9], [10], [11])

$$\frac{dH}{dt} = cH + \eta h - \varepsilon H^3$$

$$\frac{dh}{dt} = -\theta H - \gamma h \tag{8}$$

with initial conditions as in (see Ref. [8], [9], [10], [11])

$$H(0) = 1,\ h(0) = 1 \tag{9}$$

where $c, \eta, \gamma, \theta$ and $\varepsilon$ are physical constants, also $\varepsilon$ is perturbation coefficient which is defined $0 < \varepsilon < 1$ and small enough. $H$ is temperature of the eastern equatorial Pasific sea surface and $h$ is the thermo-cline depth anomaly.

Let's, $W(k)$ denotes the transformed form of $H(t)$ and $V(k)$ denotes the transformed form of $h(t)$. Also $N(k)$ denotes the transformed form of the nonlinear term $N(H(t)) = H(t)^3$. Then, transformed form of equations (8) and (9) can be written as

$$W(k+1) = \frac{cW(k) + \eta V(k) - \varepsilon N(k)}{k+1}$$

$$V(k+1) = \frac{-\theta W(k) - \gamma V(k)}{k+1} \tag{10}$$

$$W(0) = 1, V(0) = 1$$

where from (5),(6) and (7), $N(k)$ are calculated as below

$$N(0) = W^3(0)$$
$$N(1) = 3W^2(0)W(1)$$
$$N(2) = 3W(0)W^2(1) + 3W^2(0)W(2)$$
$$N(3) = 3W^2(0)W(3) + 6W(0)W(1)W(2) + W^3(1) \tag{11}$$
$$N(4) = 3W(0)W^2(2) + 3W^2(1)W(2) + 3W^2(0)W(4) + 6W(0)W(1)W(3)$$
$$N(5) = 3W(1)W^2(2) + 3W^2(1)W(3) + 3W^2(0)W(5) + 6W(0)W(1)W(4) + 6W(0)W(2)W(3)$$

Putting (11) in place to (10) and we apply the differential transform method (DTM) to the coupled nonlinear system (10), we get the some transformed terms following

$$W(0) = 1, V(0) = 1,\ W(1) = c + \eta - \varepsilon, V(1) = -\theta - \gamma$$

$$W(2) = \frac{c^2 + \eta(c - 3\varepsilon - \theta - \gamma) + \varepsilon(3\varepsilon - 4c)}{2},\ V(2) = \frac{\theta(-c - \eta + \varepsilon + \gamma) + \gamma^2}{2} \tag{12}$$

$$W(3) = \frac{c\left(c^2 + \eta c - 13c\varepsilon - 2\eta\theta - \eta\gamma - 18\varepsilon\eta + 27\varepsilon^2\right)}{6} - \frac{\eta\left(\eta\theta - 4\theta\varepsilon - \gamma\theta - \gamma^2 - 3\varepsilon\gamma - 21\varepsilon^2 + 6\varepsilon\eta\right) - 15\varepsilon^3}{6}$$

$$V(3) = \frac{\theta\left(-c^2 - \eta c + 4c\varepsilon + \eta\theta + 2\theta\gamma\right) - \gamma^3}{6} + \frac{\theta\left(3\eta\varepsilon - 3\varepsilon^2 + c\gamma - \gamma\varepsilon - \gamma^2\right)}{6}$$

⋮ ⋮

By continuing in this way and also combining (12), we obtain the DTM solutions of (8)

$$H(t) = \sum_{k=0}^{\infty} W(k)t^k \bigg|_{t_0=0} \cong 1 + (c+\eta-\varepsilon)t + \frac{(c(c-4\varepsilon+\eta)+\eta(-\theta-3\varepsilon-\gamma)+3\varepsilon^2)}{2}t^2$$

$$+ \frac{c(c^2+\eta c-13c\varepsilon-2\eta\theta-\eta\gamma-18\varepsilon\eta+27\varepsilon^2)}{6}t^3 \quad (13)$$

$$- \frac{\eta(\eta\theta-4\theta\varepsilon-\gamma\theta-\gamma^2-3\varepsilon\gamma-21\varepsilon^2+6\varepsilon\eta)-15\varepsilon^3}{6}t^3 + \cdots$$

$$h(t) = \sum_{k=0}^{\infty} V(k)t^k \bigg|_{t_0=0} \cong 1 + (-\theta-\gamma)t + \left(\frac{\theta(-c-\eta+\varepsilon+\gamma)+\gamma^2}{2}\right)t^2 \quad (14)$$

$$+ \left(\frac{\theta(-c^2-\eta c+4c\varepsilon+\eta\theta+2\theta\gamma+)-\gamma^3}{6}\right)t^3 + \frac{\theta(3\eta\varepsilon-3\varepsilon^2+c\gamma-\gamma\varepsilon-\gamma^2)}{6}t^3 + \cdots$$

which are interaction between sea surface temperature (SST) $H(t)$ and thermo-cline depth anomaly $h(t)$ as in shown the Figure 1. Also, DTM solution of some numerical values for equation (13) is obtained to verify the accuracy and effectiveness in the Table 1, Table 2 and figures 2-4.

| | DTM | | ADM | | VIM | |
|---|---|---|---|---|---|---|
| $t$ | $\varepsilon=0.1$ | $\varepsilon=0.2$ | $\varepsilon=0.1$ | $\varepsilon=0.2$ | $\varepsilon=0.1$ | $\varepsilon=0.2$ |
| 0.0 | 1.000000000 | 1.000000000 | 1.000000000 | 1.000000000 | 1.000000000 | 1.000000000 |
| 0.2 | 1.363075110 | 1.329016721 | 1.363075110 | 1.329016721 | 1.363041744 | 1.328937177 |
| 0.4 | 1.678076688 | 1.578711729 | 1.678076689 | 1.578711729 | 1.677354988 | 1.577587314 |
| 0.6 | 1.926967553 | 1.737503048 | 1.926967553 | 1.737503047 | 1.922756278 | 1.734769108 |
| 0.8 | 2.100805593 | 1.809993417 | 2.100805593 | 1.809993417 | 2.087689868 | 1.815479092 |
| 1.0 | 2.197774051 | 1.759598132 | 2.197774050 | 1.759598132 | 2.173204440 | 1.839161918 |

**Table 1:** Compare the numerical solutions of Eq. 8 for $H(t)$ obtained by DTM, ADM and VIM while $c=1$, $\eta=1, \gamma=1, \theta=1$.

| | DTM | | ADM | | VIM | |
|---|---|---|---|---|---|---|
| $t$ | $\varepsilon=0.1$ | $\varepsilon=0.2$ | $\varepsilon=0.1$ | $\varepsilon=0.2$ | $\varepsilon=0.1$ | $\varepsilon=0.2$ |
| 0.0 | 1.000000000 | 1.000000000 | 1.000000000 | 1.000000000 | 1.000000000 | 1.000000000 |
| 0.2 | 1.618047013 | 1.568552948 | 1.618047013 | 1.568552948 | 1.618050569 | 1.568813651 |
| 0.4 | 2.287634977 | 2.094478815 | 2.287634977 | 2.094478814 | 2.291096651 | 2.100926397 |
| 0.6 | 2.922322488 | 2.488466817 | 2.922322488 | 2.488466817 | 2.953611643 | 2.526631543 |
| 0.8 | 3.391665078 | 2.950478413 | 3.391665078 | 2.950478416 | 3.534733325 | 2.881459315 |
| 1.0 | 3.258944085 | 7.752628428 | 3.258944090 | 7.752628479 | 4.030556145 | 3.272412003 |

**Table 2:** Compare the numerical solutions of Eq. 8 for $H(t)$ obtained by DTM, ADM and VIM while $c=2$, $\eta=1, \gamma=1, \theta=1$.

**Problem 2:** Secondly, by making some changes on the equation (8), we consider the following Enso delayed oscillator model (see more Ref. [12],[13],[14]);

$$(1-\beta\sigma)\left(\frac{dH}{dt}\right) = (\alpha-\beta)H - \varepsilon H^3 \quad (15)$$

where $\sigma$ is time delayed and $\eta = 0$ in (8). This equation is the oscillator model which emphasizes the ocean and atmosphere interactions in the equatorial eastern Pacific and anomaly variations only in this coupling region ( see Ref. [12],[13],[14]). In eq. (15), $\sigma$, $\beta$, $\alpha$ and $\varepsilon$ are positive physical variables. Again, we apply the procedures (10) and (12) to equation (15). Thus, we obtain the differential transform of (15) and some transformed terms of (15) as follows respectively

$$W(k+1) = \frac{(\alpha-\beta)}{(1-\beta\sigma)(k+1)}W(k) - \frac{\varepsilon}{(1-\beta\sigma)(k+1)}N(k) \quad (16)$$

$$W(0) = 1$$

and

$$W(0)=1, W(1)=-\frac{\alpha-\beta-\varepsilon}{\beta\sigma-1}, W(2)=\frac{1}{2}\frac{(\alpha-\beta-\varepsilon)(\alpha-\beta-3\varepsilon)}{(\beta\sigma-1)^2}$$

$$W(3) = -\frac{1}{6}\frac{(\alpha-\beta-\varepsilon)(\alpha^2-2\alpha\beta-12\alpha\varepsilon+\beta^2+12\beta\varepsilon+15\varepsilon^3)}{(\beta\sigma-1)^3} \quad (17)$$

$$W(4) = \frac{1}{24} \frac{(\alpha-\beta-\varepsilon)(\alpha^3 - 3\alpha^2\beta - 39\alpha^2\varepsilon + 3\alpha\beta^2 + 78\alpha\beta\varepsilon + 135\alpha\varepsilon^2)}{(\beta\sigma-1)^4}$$

$$+ \frac{1}{24} \frac{(\alpha-\beta-\varepsilon)(-\beta^3 - 39\beta^2\varepsilon - 135\beta\varepsilon^2 - 105\varepsilon^3)}{(\beta\sigma-1)^4}$$

$$\vdots \qquad \vdots \qquad \vdots$$

By pursuing this process and also unifying (17), we obtain the DTM solutions of (15) as follow

$$H(t) = \sum_{k=0}^{\infty} W(k)t^k \bigg|_{t_0=0} \cong 1 - \frac{\alpha-\beta-\varepsilon}{\beta\sigma-1}t + \frac{1}{2}\frac{(\alpha-\beta-\varepsilon)(\alpha-\beta-3\varepsilon)}{(\beta\sigma-1)^2}t^2$$

$$- \frac{1}{6}\frac{(\alpha-\beta-\varepsilon)(\alpha^2 - 2\alpha\beta - 12\alpha\varepsilon + \beta^2 + 12\beta\varepsilon + 15\varepsilon^3)}{(\beta\sigma-1)^3}t^3 \qquad (18)$$

$$+ \frac{t^4}{24(\beta\sigma-1)^4}(\alpha-\beta-\varepsilon)(\alpha^3 - 3\alpha^2\beta - 39\alpha^2\varepsilon + 3\alpha\beta^2 + 78\alpha\beta\varepsilon + 135\alpha\varepsilon^2)$$

$$+ \frac{1}{24(\beta\sigma-1)^4}(\alpha-\beta-\varepsilon)(-\beta^3 - 39\beta^2\varepsilon - 135\beta\varepsilon^2 - 105\varepsilon^3)t^4 + \cdots$$

By taking numerical values for $\alpha, \beta, \theta, \sigma$ and $\varepsilon$, compared results for equation (18) is given to confirm the accuracy and effectiveness of DTM solutions in Table 3, Table 4 and in figure 5 to figure 10 as below.

|  | Exact | | DTM | | ADM | | VIM | |
|---|---|---|---|---|---|---|---|---|
| $t$ | $\varepsilon=0.05$ | $\varepsilon=0.1$ | $\varepsilon=0.05$ | $\varepsilon=0.1$ | $\varepsilon=0.05$ | $\varepsilon=0.1$ | $\varepsilon=0.05$ | $\varepsilon=0.1$ |
| 0.0 | 1.000000000 | 1.000000000 | 1.000000000 | 1.000000000 | 1.000000000 | 1.000000000 | 1.000000000 | 1.000000000 |
| 0.4 | 1.065476869 | 1.042243490 | 1.065476869 | 1.042243489 | 1.065476869 | 1.042243491 | 1.065440029 | 1.042238926 |
| 0.8 | 1.131796087 | 1.082251484 | 1.131796089 | 1.082251485 | 1.131796088 | 1.082251485 | 1.131713720 | 1.082260100 |
| 1.2 | 1.198337436 | 1.119727591 | 1.198337435 | 1.119727591 | 1.198337435 | 1.119727591 | 1.198221187 | 1.119745754 |
| 1.6 | 1.264434686 | 1.154458882 | 1.264434686 | 1.154458882 | 1.264434686 | 1.154458882 | 1.264317518 | 1.154460568 |
| 2.0 | 1.329402528 | 1.186318727 | 1.329402528 | 1.186318729 | 1.329402529 | 1.186318730 | 1.329331004 | 1.186269275 |

**Table 3:** Compare the numerical solutions of Eq. 15 obtained by DTM, ADM and VIM while $\alpha=0.5, \beta=0.3, \sigma=0.25$.

| | Exact | | DTM | | ADM | | VIM | |
|---|---|---|---|---|---|---|---|---|
| $t$ | $\varepsilon=0.05$ | $\varepsilon=0.1$ | $\varepsilon=0.05$ | $\varepsilon=0.1$ | $\varepsilon=0.05$ | $\varepsilon=0.1$ | $\varepsilon=0.05$ | $\varepsilon=0.1$ |
| 0.0 | 1.000000000 | 1.000000000 | 1.000000000 | 1.000000000 | 1.000000000 | 1.000000000 | 1.000000000 | 1.000000000 |
| 0.4 | 1.261904652 | 1.222317642 | 1.261904653 | 1.222317640 | 1.261904652 | 1.222317640 | 1.254384957 | 1.217730500 |
| 0.8 | 1.562236074 | 1.450459383 | 1.562236064 | 1.450459413 | 1.562236064 | 1.450459413 | 1.540458840 | 1.442069517 |
| 1.2 | 1.884088435 | 1.663165659 | 1.884085615 | 1.663174675 | 1.884085615 | 1.663174675 | 1.847702160 | 1.657228153 |
| 1.6 | 2.200074136 | 1.841919059 | 2.199919220 | 1.842415020 | 2.199919218 | 1.842415017 | 2.160603481 | 1.845724524 |
| 2.0 | 2.480426774 | 1.978020657 | 2.477167800 | 1.988570079 | 2.477167800 | 1.988570063 | 2.459058672 | 1.991364765 |

**Table 4:** Compare the numerical solutions of Eq. 15 obtained by DTM, ADM and VIM while $\alpha=1, \beta=0.5, \sigma=0.5$.

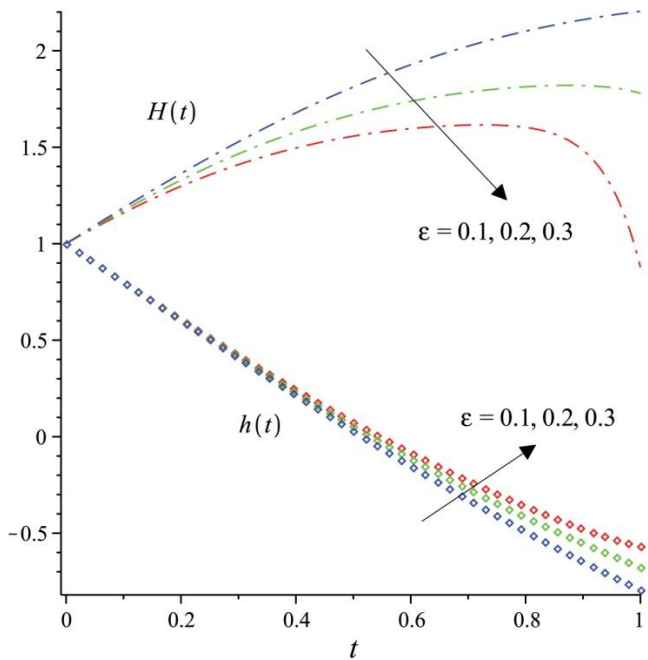

**Figure 1:** Effects of different values of $\varepsilon$ on the interaction between $H(t)$ and $h(t)$ for equation

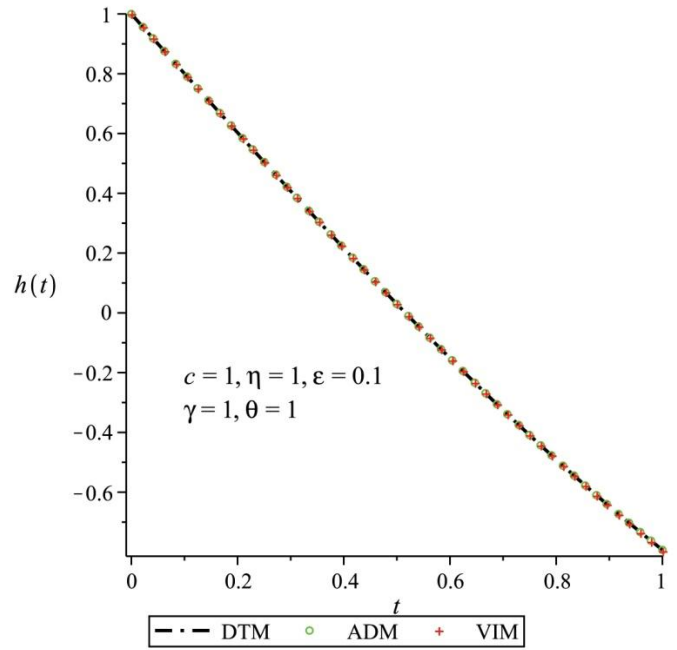

**Figure 2:** Comparison for the solutions of thermo-cline depth anomaly $h(t)$ for equation (8).

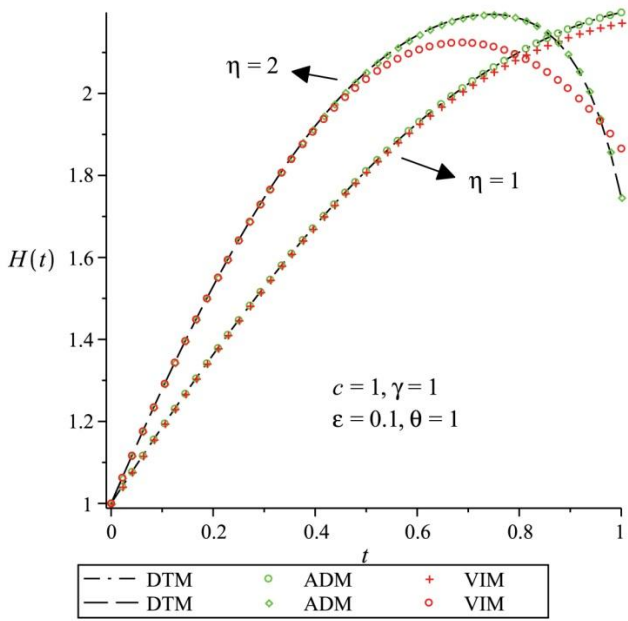

**Figure 3**: Comparison of DTM, ADM and VIM solutions for $H(t)$ in eq. (8) and effects of different values of $\eta$.

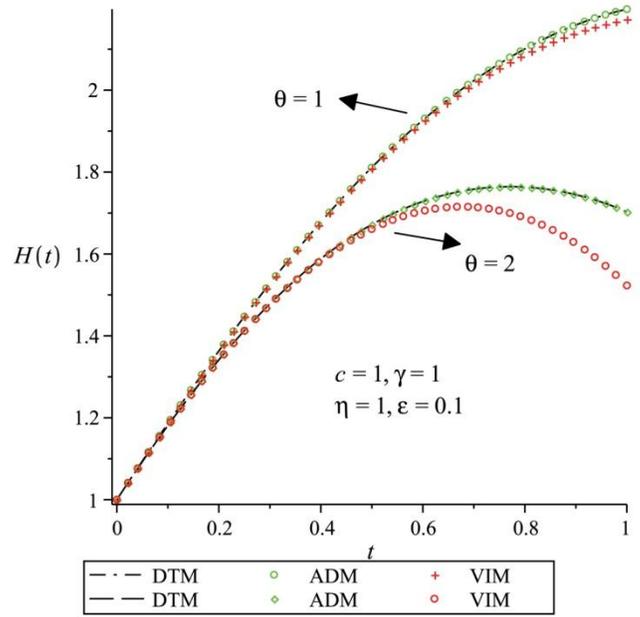

**Figure 4:** Comparison of DTM, ADM and VIM solutions for $H(t)$ in eq. (8) and effects of different values of $\theta$.

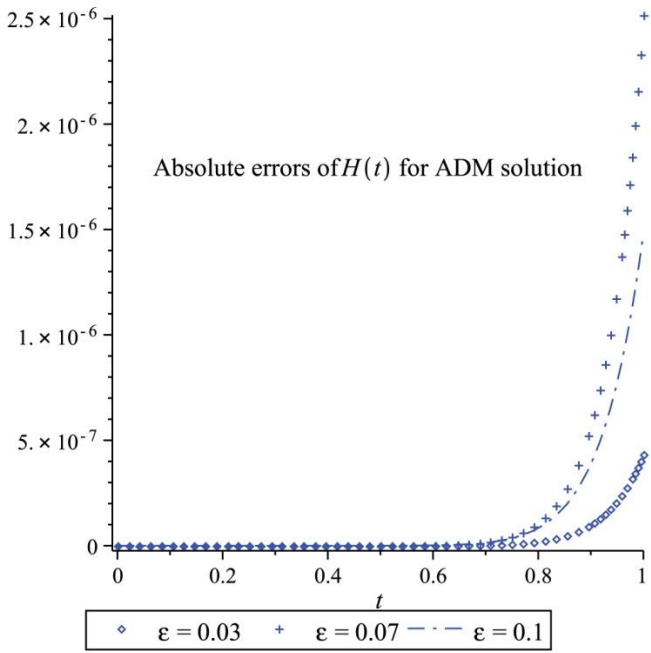

**Figure 5:** Error values of eq. (15) obtained by ADM for different values of $\varepsilon$.

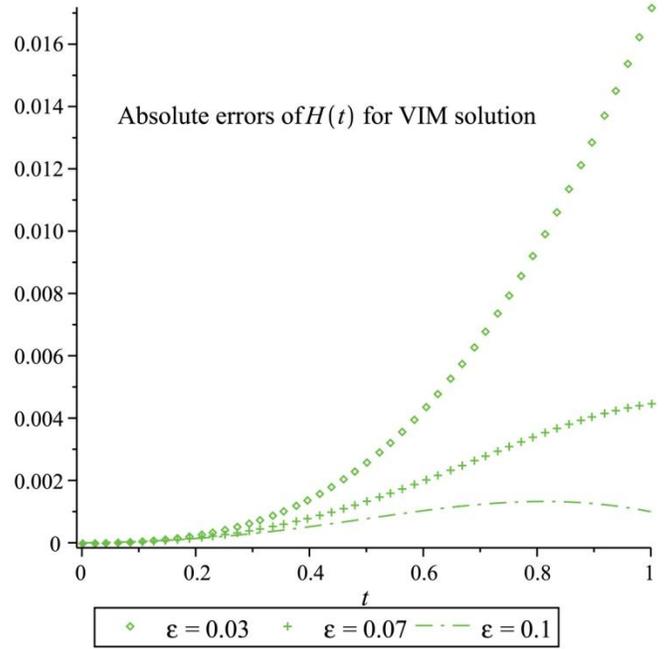

**Figure 6:** Error values of eq. (15) obtained by VIM for different values of $\varepsilon$.

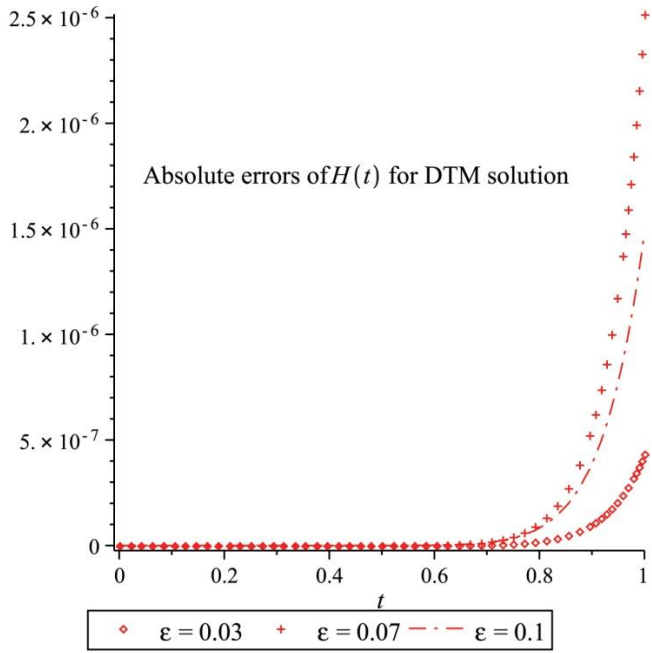

**Figure 7:** Error values of eq. (15) obtained by DTM for different values of $\varepsilon$.

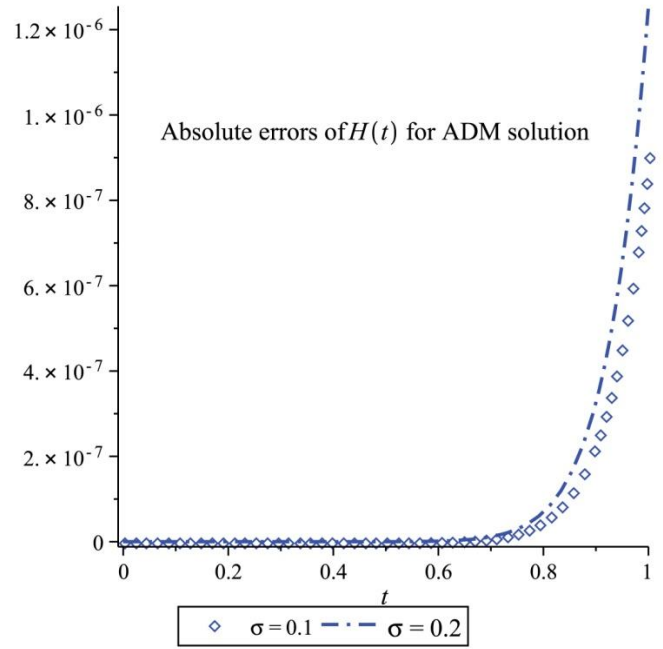

**Figure 8:** Error values of eq. (15) obtained by ADM for different values of $\sigma$.

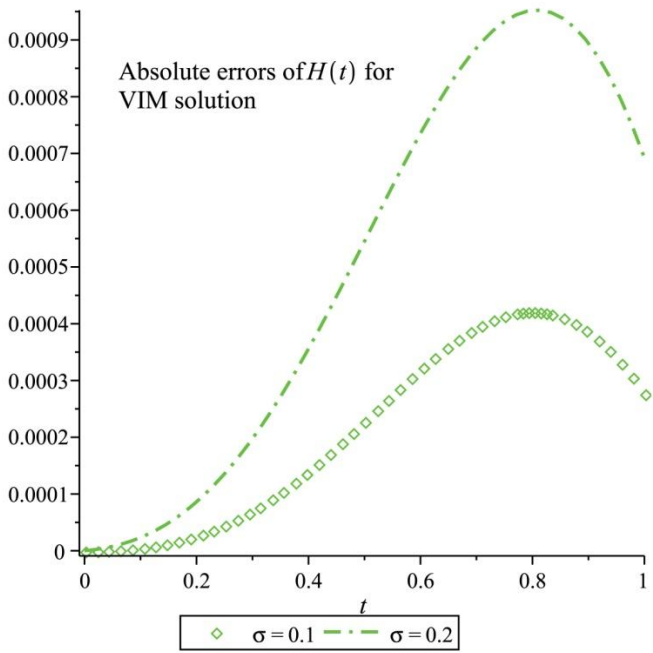

**Figure 9:** Error values of eq. (15) obtained by VIM for different values of $\sigma$.

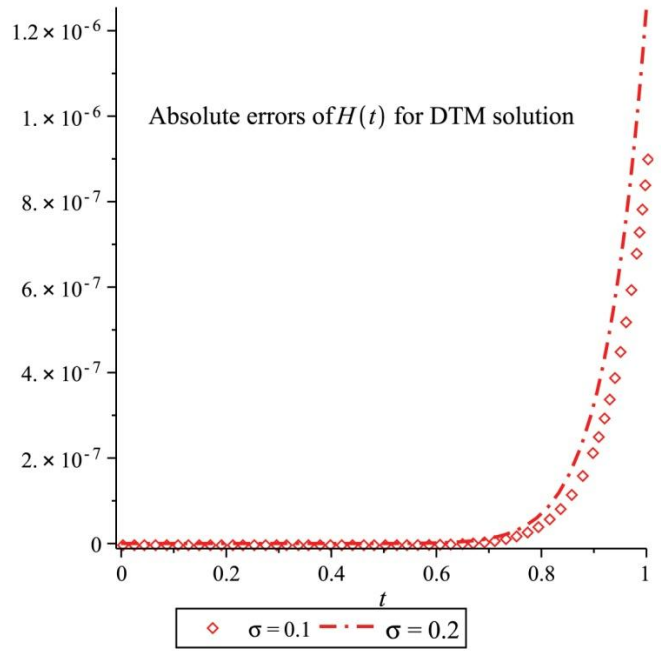

**Figure 10:** Error values of eq. (15) obtained by DTM for different values of $\sigma$.

## 4. Conclusions

Differential transform method (DTM) has been successfully applied to the two ENSO nonlinear model. DTM, VIM and ADM solutions are compared. It's noted that we have obtained fairly good results by DTM and ADM with compare to VIM as shown in Table 1-4. Accuracy and errors are demonstrated in Figure 2 to 10 that DTM is powerful semi-analytic method to solve nonlinear ordinary differential system.

We consider Eq. 8 in terms of interaction between sea surface temperature and thermo-cline depth anomaly in Figure 1. In Figure 2-4, thermo-cline depth anomaly and sea surface temperature are investigated separately with respect to various values of $\theta, \eta$. Additionally, we obtain errors of solution for Eq. 15 according to the changing values of $\varepsilon, \sigma$ in Figure 5 to 10.